\documentclass[12pt]{amsart}

\usepackage{verbatim}
\usepackage{amssymb}
\usepackage{amsfonts}
\usepackage{amsbsy}

\newcommand{\esubproof}{\hfill$\square$}
\newcommand{\Proof}{\noindent {\bf Proof.~~}}
\newcommand{\R}{{\mathbb R}}

\newcommand{\N}{{\mathbb N}}

\newtheorem{prop}{Proposition}[section]
\newtheorem{lem}[prop]{Lemma}

\newtheorem{claim}[prop]{Claim}

\newtheorem{coro}[prop]{Corollary}
\newtheorem{theo}[prop]{Theorem}

\begin{document}

\title[Arcwise connectedness]
{Arcwise connectedness of the boundaries of connected self-similar
sets}

\author[T. -M. Tang]{Tai-Man Tang}
\address{Department of Mathematics\\ The Chinese University of Hong Kong\\ Shatin,
NT \\ Hong Kong.}
\thanks{This research is supported by Prof. Ka-Sing Lau.}
\email{tmtang@math.cuhk.edu.hk}

\begin{abstract}
Let $T$ be the attractor of injective contractions
$f_1,\ldots,f_m$ on $\R^2$ that satisfy the Open Set Condition. If
$T$ is connected, $\partial T$ is arcwise connected. In
particular, the boundary of the L\'{e}vy dragon is arcwise
connected.

\vspace{1ex} \noindent {\bf Key Words:} self-similar sets,
L\'{e}vy dragon, L\'{e}vy curve, reptiles, self-affine tiles.

\vspace{1ex} \noindent {\bf AMS subject classification (2000):}
 28A80(Primary); 54F65(Secondary)

\end{abstract}

\maketitle

\section{The theorem}

Let $f_1,\ldots,f_m$ be a family of injective contractions on
$\R^2$ satisfying the Open Set Condition: there is a nonempty
bounded open set $V$ such that $f_i(V)\cap f_j(V)=\emptyset$ for
$i\neq j$, and $\cup_{i=1}^mf_i(V)\subset V$ (see e.g. \cite{F}).
Let $T$ be the attractor of the system. Suppose that $T$ is
connected. Among other results, Luo, Rao and Tan prove that
$\partial T$ is connected \cite[Theorem 1.1]{LRT}. They further
ask whether $\partial T$ is arcwise connected. We answer the
question in the affirmative.

\vspace{5mm}

\begin{theo}\label{theo1.1} Let $f_1,\ldots, f_m$ be a family of
injective contractions on $\R^2$ satisfying the Open Set
Condition. Suppose that $T$ is connected. Then $\partial T$ is
arcwise connected.
\end{theo}

\begin{coro}\label{coro1.2} The boundary of a connected reptile or self-affine tile
is arcwise connected.
\end{coro}

Lately there is some interest in the topology of self-similar
sets, particularly for some classical reptiles and self-affine
tiles (see \cite{BKS}, \cite{BW}, \cite{LRT}, \cite{NN}). If the
$f_i$ are similarities of the same contraction ratio and $T^\circ
\neq \emptyset$, $T$ is called a reptile. A self-affine tile is
defined by an expanding matrix and a digit set. The twindragon,
the Heighway dragon and the L\'{e}vy dragon are classical examples
in both classes.  Bandt and Wang \cite{BW} show that the
twindragon is a disk.  Ngai and Nguyen \cite{NN} show that the
Heighway dragon is a union of disks, each having a common point
with each of its two neighboring disks. Hence our theorem is true
for the twindragon and the Heighway dragon. Notice that for these
$T$, almost all points in $\partial T$ are boundary points of the
components of $T^\circ$. The only exceptions are the two special
points of the Heighway dragon, which are limit points of such
components.

The non-trivial cases for our theorem are offered by those $T$
where $\partial T$ has many points that are not boundary points of
the components of $T^\circ$, but are the limits of such
components. The L\'{e}vy dragon offers an example. The Hausdorff
dimension of its boundary has been calculated using different
methods \cite{DK}, \cite{SW}. Its topology is discussed by Bailey,
Kim and Strichartz \cite{BKS}. The arcwise connectedness of its
boundary is an addition to the results there.

\section{Preliminaries}

We collect here some definitions and results from point set
topology and self-similar sets.

A {\em continuum} is a compact connected set. It is {\em
non-degenerate} if it has more than one point. Let $S$ be a
topological space. Let $G$ be an infinite collection of subsets of
$S$, not necessarily different. The set of $x \in S$ such that
every neighborhood of $x$ contains points of infinitely many sets
in $G$ is called the {\em limit superior} of $G$, denoted $\limsup
G$. The set of $y \in S$ such that every neighborhood of $y$
contains points from all but a finite number of the sets of $G$ is
called the {\em limit inferior} of $G$, written $\liminf G$. If
$\liminf G=\limsup G$, then $G$ is said to be convergent, with
limit $\lim G = \liminf G=\limsup G$.

A set $M$ is said to be {\em locally connected} at $p\in M$ if for
every neighborhood $U$ of $p$, there exists a neighborhood $V$ of
$p$ such that every point of $M\cap V$ lies in the component of
$M\cap U$ containing $p$. Equivalently, $M$ has a local base at
$p$ consisting of connected sets. $M$ is {\em locally connected}
if it is locally connected at every one of its points.

\begin{theo}\label{theo2.1} (a) \cite[p.13]{W} If the continuum $M\subset \R^2$ is
not locally connected at one of its points $p$, then there is a
ball $B_r(p)$ and an infinite sequence of distinct components $C$
and $C_i$ of $M\cap \overline{B_r(p)}$, $i=1,2,\ldots$, such that
$\lim \{C_i\}=C$ and $p\in C$.

(b) \cite[p.14]{W} There is a non-degenerate subcontinuum $H$ of
$M$ containing $p$ such that $M$ is not locally connected at every
point of $H$.
\end{theo}

\begin{theo}\label{theo2.2} \cite[p.27]{W} Every locally connected continuum
is arcwise connected.
\end{theo}

An {\em arc} is a homeomorphic image of $[0,1]$. A {\em simple
closed curve} is a homeomorphic image of a circle. A set $M$ is
said to have {\em property S} if for each $\epsilon>0$, $M$ is the
union of a finite number of connected sets of diameter less than
$\epsilon$.

\begin{theo}\label{theo2.3} \cite[p.19]{W} A continuum $M$ is locally
connected if and only if $M$ has property $S$.
\end{theo}

\begin{theo}\label{theo2.4}  \cite[p.34]{W}  If $M \subset
\R^2$ is a locally connected continuum with no cut point, the
boundary of any component of $\R^2 \setminus M$ is a simple closed
curve.
\end{theo}

For $\alpha = i_1\ldots i_k \in \{1,\ldots,m\}^k$,  we write
$f_\alpha=f_{i_1}\circ \cdots \circ f_{i_k}$. $f_\alpha (T)$ is
called a $k$th-level piece of $T$. Let
\[  c_i=\sup \{ \frac{|f_i(x)-f_i(y)|}{|x-y|}: x\neq y;
x,y \in \R^2\}<1; c=\sup\{c_1,\ldots,c_m\}<1. \] Notice that
$\mbox{diam}(f_\alpha(T))\leq c^k\mbox{diam}(T)\rightarrow 0$ as
$k\rightarrow \infty$. As $T=\cup_{\alpha \in
\{1,\ldots,m\}^N}f_\alpha(T)$, and each $f_\alpha(T)$ is
connected, we have part (a) of the following.

\begin{theo}\label{theo2.5} Let $T$ be the connected attractor of injective
contractions $f_i$, as in Theorem \ref{theo1.1}. Then

(a) $T$ has property $S$.

(b) $T$ is arcwise connected (\cite{H}, \cite[p.33]{K}).

(c) {\em If $T^\circ \neq \emptyset$, $T=\overline{T^\circ}$}
(e.g. \cite[p.226]{LRT}).

(d) {\em Suppose $T^\circ \neq \emptyset$. For different
$\alpha_1, \alpha_2 \in \{1,\ldots m\}^k$, $f_\alpha(T^\circ)\cap
f_\beta(T^\circ)=\emptyset$} (e.g. \cite[p.226]{LRT}).
\end{theo}

\section{The Proof}

We prove Theorem \ref{theo1.1} in this section. Under the given
hypothesis, $\partial T$ is connected (\cite[Theorem
1.1(ii)]{LRT}) and hence a continuum. We will prove that it is
arcwise connected.

\begin{lem}\label{lem3.1} If $T^\circ=\emptyset$, $\partial T$ is arcwise connected.
\end{lem}

\Proof In this case $\partial T=T$. The arcwise connectedness of
$\partial T$ follows from that of $T$ (Theorem \ref{theo2.5}(b)).
\esubproof

Hereafter, we assume that $T^\circ \neq \emptyset$. Suppose
$\partial T$ is not arcwise connected. We derive a contradiction
in a sequence of steps.

\begin{claim}\label{claim3.2} Suppose that $\partial T$ is not arcwise connected.
There is a point $p\in \partial T$, and an open ball $B_r(p)$ such
that $\partial T \cap \overline{B_r(p)}$ has infinitely many
components $C$ and $C_i$, $i=1,2,\ldots$, such that
$\lim\{C_i\}=C$ and $p \in C$.
\end{claim}

\Proof $\partial T$ not arcwise connected implies that it is not
locally connected (Theorem \ref{theo2.2}). The result follows from
Theorem \ref{theo2.1}(a).
\esubproof

\begin{claim}\label{claim3.3} Let $N$ be a positive integer such that for any $N$th-level
piece $f_\alpha(T)$, $\alpha \in \{1,\ldots,m\}^N$,
$\mbox{diam}(f_\alpha(T))< r/2$. There is an $N$th-level piece of
$T$, denoted $A$, that is contained in $B_r(p)$ and intersects
infinitely many $C_i$.
\end{claim}

\Proof  As $\lim \{C_i\} = C$, $C_i \cap B_{r/2}(p)\neq \emptyset$
except for finitely many $i$. As $C_i \subset
\partial T \subset T=\cup_{\alpha\in
\{1,\ldots,m\}^N}f_\alpha(T)$, each of these points of
intersections is in some $N$th-level piece of $T$. As only
finitely many of such pieces intersect $B_{r/2}(p)$, some
$N$th-level piece of $T$, called $A$, contains points from
$B_{r/2}(p)\cap C_i$ for infinitely many $i$. As
$\mbox{diam}(A)<r/2$, $A \subset B_r(p)$.
\esubproof

Blow up $T$ so that every point in $\partial T$ is an interior
point of the blow up. In detail, choose a $k$th-level piece
$f_\alpha(T)$ of $T$ with $\partial (f_\alpha(T))\subset T^\circ$.
The blow up is $f_\alpha^{-1}(T)$, in which $T$ is bordered by
neighbors, $f_\alpha^{-1}$ of the other $k$th-level pieces of $T$.

\begin{claim}\label{claim3.4} There is an $N$th-level piece of a neighbor of $T$, called
$B$, such that $B \subset B_r(p)$ and $B\cap A$ contains points
from two of the $C_i$'s, say $C_1, C_2$. Here $A$ is as in Claim
\ref{claim3.3}
\end{claim}

\Proof Choosing another $N$ if necessary, suppose that the
$N$th-level pieces of $T$ and its neighbors in the blow up have
diameter less than $r/2$. From Claim \ref{claim3.3}, $A \cap C_i
\cap B_{r/2}(p)\neq \emptyset$ for infinitely many $i$. As
$C_i\subset
\partial T$, $A \cap C_i \cap B_{r/2}(p)$ is also contained in the
neighbors of $T$. As only finitely many $N$th-level pieces of the
neighbors of $T$ intersects $B_{r/2}(p)$, one such piece $B$
contains points in $A\cap C_i\cap B_{r/2}(p)$ for infinitely many
$i$. As $\mbox{diam}(B)<r/2$, $B \subset B_r(p)$.

By renaming the $C_i$'s if necessary, suppose that $A \cap B$
contains points from $C_1, C_2$. \esubproof

Let $x \in A\cap B \cap C_1$, $y \in A\cap B\cap C_2$. As $A$ and
$B$ are arcwise connected (Theorem \ref{theo2.5}(b)), there are
arcs $\gamma \subset A \subset T$, $\beta \subset B \subset
\overline{T^c}$ with endpoints $x,y$. We get a contradiction by
proving the following.

\begin{claim}\label{claim3.5} $C_1$ and $C_2$ cannot be distinct components of $\partial T
\cap \overline{B_r(p)}$.
\end{claim}

\Proof {\it Case 1.} If $\gamma = \beta$, the arcs are in
$\partial T \cap \overline{B_r(p)}$, and the claim is true.

{\it Case 2.} Suppose that $\gamma \neq \beta$, and $\gamma \cap
\beta=\{x,y\}$. That is, $\gamma \cup \beta$ is a simple closed
curve enclosing a region $D \subset B_r(p)$.

If $\gamma \subset \partial T$ or $\beta\subset \partial T$, then
$C_1$ and $C_2$ are joined by an arc in $\partial T \cap
\overline{B_r(p)}$, and the claim is true.

Suppose that $\gamma \cap T^\circ \neq \emptyset$. Look at the
components of $T^\circ \cap D$ whose boundary has nonempty
intersection with $\gamma$. Call them $A_i$, $i \in \N$. Notice
that $\partial A_i \subset (\partial T \cap D) \cup \gamma \cup
\beta$.

We claim that $\overline{A_i}$ is a locally connected continuum
with no cut point. We have to prove the local connectedness of
$\overline{A_i}$  at each of its points. As $T$ is a locally
connected continuum (Theorem \ref{theo2.3}, \ref{theo2.5}(a)), it
is locally connected at each of its points. For $z\in D\cap
\overline{A_i}$, local connectedness of $\overline{A_i}$ at $z$
follows from the local connectedness of $T$ at $z$.

Next consider $z \in \overline{A_i}\cap \gamma$ with the property
that there is an interval $(t_1,t_2)\subset [0,1]$ with $z \in
\gamma(t_1,t_2)\subset \overline{A_i}$ (the `interior boundary
points'). We have used the same symbol for the arc $\gamma$ and
one of its parametrizations $\gamma:[0,1]\rightarrow \R^2$. Notice
that $\mbox{dist}(z,\partial T \cap D)>0$. Suppose that
$\overline{A_i}$ is not locally connected at $z$. Then there is a
closed ball $S$ of $z$, such that $S\cap
\partial A_i \subset \gamma$, and $\overline{A_i}\cap S$ has
components $C_i',C'$ such that $\lim \{C_i'\}=C'$ (Theorem
\ref{theo2.1}(a)). By our choice of $S$, $\partial C_i' \subset
\gamma$. It follows that every neighborhood of $z$ in $S$
intersects $\gamma$ in infinitely many components. Hence $\gamma$
does not have a local base of connected neighborhoods at $z$,
contradictory to the local connectedness of $\gamma$. Hence
$\overline{A_i}$ is locally connected at $z$. The same argument
apply to the `interior boundary points' on
$\overline{A_i}\cap\beta$.

It remains to establish the locally connectedness of
$\overline{A_i}$ at the `corner boundary points' of
$\overline{A_i}$, the points $z = \gamma(t) \in \gamma$ (and the
similar points on $\beta$) with the following property. There is
no interval $(t_1,t_2) \subset [0,1]$ containing $t$ such that
$\gamma(t_1,t_2) \subset \overline{A_i}$. If $\overline{A_i}$ is
not locally connected at $z$, it is not locally connected on a
non-degenerate sub-continuum $H$ of $\overline{A_i}$ containing
$z$ (Theorem A(b)). As we have established the local connectedness
of $\overline{A_i}$ at the points of $\overline{A_i}$ in $D$, $H
\subset \gamma$ and hence must be a non-degenerate sub-arc. But
then points in $H$ other then its two end points are the `interior
boundary points' discussed in the last paragraph, and
$\overline{A_i}$ is locally connected at such points. This
contradicts the definition of $H$, and proves the local
connectedness of $\overline{A_i}$ at $z$. Hence $\overline{A_i}$
is locally connected.

$\overline{A_i}$ has no cut point, as for any $z \in
\overline{A_i}$, $A_i\setminus \{z\}$ is in one component, and
hence so is $\overline{A_i}\setminus \{z\}$. This establishes our
claim that $\overline{A_i}$ is a locally connected continuum with
no cut point.

By Theorem \ref{theo2.4}, the boundaries of the components of
$\R^2 \setminus \overline{A_i}$ are simple closed curves.  Let
$\delta_i$ be the boundary of the unbounded component. Points on
$\delta_i$ are of three types: those in $D$, $\gamma$ or $\beta$.
Those in $D$ and $\beta$ are in $\partial T$.

Let $s_i:= \inf \{s: \gamma(s) \in \delta_i\}$, $t_i:=\sup \{ t:
\gamma(t) \in \delta_i\}$. Then $\delta_i\setminus \{\gamma(s_i),
\gamma(t_i)\}$ is consist of two parts, with at least one lying
entirely in $D \cup \beta$. Call one such part $\delta_i'$. Then
$\delta_i'\subset  \partial T$. Define \[  \gamma':=(\gamma
\setminus \bigcup_i\gamma(s_i,t_i)) \cup (\bigcup_i \delta_i'). \]
Then $\gamma' \subset \partial T$. Though $\gamma'$ may not be an
arc, it is the image of a continuous curve joining $x,y$.
Therefore $x, y$ and hence $C_1, C_2$ are in the same component of
$\partial T \cap \overline{B_r(p)}$. This finishes the argument
when $\gamma \cap \beta = \{x,y\}$.

{\it Case 3.} Suppose that $\gamma \neq \beta$ and $\gamma \cap
\beta$ is more than $\{x,y\}$. Let $(u_i,v_i)$, $i\in \N$, be
maximal intervals with $\gamma(u_i,v_i)\cap \beta = \emptyset$.
For each $i$, $\gamma(u_i), \gamma(v_i)$ bounds a segment from
each of $\gamma$ and $\beta$. The two segments bounded a region
$D_i$. Apply the argument in case 2 to get a curve in $\partial T
\cap \overline{D_i}$ joining $\gamma(u_i)$ and $\gamma(v_i)$.
Together with the observation that $\gamma \cap \beta \subset
\partial T \cap \overline{B_r(p)}$, we get that $x,y$ and
$C_1, C_2$ are in the same component of $\partial T \cap
\overline{B_r(p)}$. \esubproof

The contradication obtained in Claim \ref{claim3.5} proves the
arcwise connectedness of $\partial T$.

\end{document}